%01/07/00
\magnification 1200
\input amstex
\documentstyle{amsppt}
\voffset -0.53cm
\vsize=7.95in
\topmatter
\rightheadtext{Sur l'application de sym\'etrie CR formelle}
\leftheadtext{Jo\"el Merker}
\title
\'Etude de la r\'egularit\'e analytique de l'application de sym\'etrie
CR formelle
\endtitle
\author Jo\"el Merker
\endauthor
\address 
Laboratoire d'Analyse, Topologie et Probabilit\'es,
Centre de Math\'ematiques et d'Informatique, UMR 6632, 39 rue Joliot
Curie, F-13453 Marseille Cedex 13, France
Tel: 00 33 (0)4 91 11 35 50 \ \
Fax: 00 33 (0)4 91 11 35 52
\endaddress
\email 
merker\@cmi.univ-mrs.fr 00 33 / (0)4 91 11 36 72 / (0)4 91 53 99 05
\endemail
\subjclass
(revised 2000).
32H02, 32C16, 32V40, 32V35,
\endsubjclass

\def\N{{\Bbb N}}

\def\R{{\Bbb R}}
\def\C{{\Bbb C}}
\define \dl{[\![}
\define \dr{]\!]}
\def\dim{\hbox{dim}}
\def\codim{\hbox{codim}}
\def\sumg{\sum_{\gamma\in \N^m}}
\def\sumstg{\sum_{\gamma\in \N_*^m}}
\def\1{{\text{\bf 1}}}
\def\v{\vert}
\def\n{\vert\vert}

\endtopmatter

\document

\head 
\S1. Introduction
\endhead

La recherche de formes normales ({\it cf.}~travaux de Chern-Moser [7],
Moser-Webster [17], Webster [20], Huang-Krantz [13], Gong [12] et
Ebenfelt [11]) pour les sous-vari\'et\'es analytiques r\'eelles de
$\C^n$ soul\`eve la question de la convergence des normalisations
formelles. Moser et Webster ont donn\'e des exemples de surfaces
${\Cal C}^\omega$ dans $\C^2$ \`a {\it tangente complexe
isol\'ee} et {\it hyperboliques} au sens de Bishop, formellement mais
non holomorphiquement normalisables (\`a cause d'un ph\'enom\`ene de
petits diviseurs, {\it voir} [17]), m\^eme lorsque la forme normale
est elle-m\^eme analytique ou alg\'ebrique. En revanche, il
appara\^{\i}t qu'un tel ph\'enom\`ene ne se produit pas pour les
objets CR, d'apr\`es les r\'esultats r\'ecents de Baouendi-Rothschild
obtenus en collaboration avec Ebenfelt ou avec Zaitsev, et \'enonc\'es
avec des hypoth\`eses g\'en\'eriques de non-d\'eg\'en\'erescence ({\it
voir} [2,3,4,6]). Ces auteurs \'etablissent notamment qu'une
application formelle inversible entre deux sous-vari\'et\'es
CR-g\'en\'eriques $\Cal C^\omega$ finiment non-d\'eg\'en\'er\'ees (ou
essentiellement finies) et {\it minimales} de $\C^n$ est
convergente. On d\'emontre ici un th\'eor\`eme de convergence plus
g\'en\'eral, valable sans hypoth\`ese de non-d\'eg\'en\'erescence, et
qui confirme la rigidit\'e du cas CR. Ce r\'esultat s'interpr\`ete
alors comme un principe de sym\'etrie de Schwarz formel pour les
applications CR et il borne aussi le degr\'e de transcendance
de l'application formelle par rapport au corps des
fractions rationnelles ({\it cf.} Coupet-Pinchuk-Sukhov [8,9]).

\head \S2. Application de sym\'etrie et \'equivalences formelles de
vari\'et\'es CR\endhead

\subhead 2.1. Application de sym\'etrie \endsubhead
Soit $h: (M,p)\to_{\Cal F}(M',p')$ une \'equivalence formelle entre
deux sous-vari\'et\'es analytiques r\'eelles ($\Cal C^\omega$)
CR-g\'en\'eriques de $\C^n$. On suppose $p=p'=0$ dans des
coordonn\'ees centr\'ees $t\in\C^n$, $t'\in\C^n$ et on note
$m:=\dim_{CR} M=m'$, $d:=\codim_\R M=d'$, $m+d=n$. Dans un travail
ant\'erieur, l'auteur a remarqu\'e l'existence d'un invariant plus
g\'en\'eral que $h$, qu'il convient d'appeler {\it application de
sym\'etrie} et sugg\'er\'e l'int\'er\^et d'\'etablir sa
r\'egularit\'e, lorsque $(M',p')$ est holomorphiquement
d\'eg\'en\'er\'ee. Cette application {\it synth\'etise en une seule
s\'erie formelle le ``jet d'ordre infini en $\bar\nu'$ de la
vari\'et\'e de Segre complexifi\'ee conjugu\'ee} ({\it voir} [13])
$\underline{\Cal S}_{h(t)}'$'': $(\bar\nu',t)\mapsto
j_{{\bar\nu'}}^\infty \underline{\Cal S}_{h(t)}'$, $\bar\nu'\in \C^n$,
$\bar\nu'\in \underline{\Cal S}_{h(t)}'$, de la mani\`ere suivante. En
coordonn\'ees locales $t'=(w',z')\in
\C^m\times \C^d$, $\tau'=(\zeta',\xi')$, telles que les $d$
\'equations de la complexifi\'ee extrins\`eque $(\Cal M',0):=
((M')^c,0^c)$ s'\'ecrivent $\xi'=\Theta'(\zeta',t')=\sum_{\gamma\in
\N^m} {\zeta'}^{\gamma} \, \Theta_\gamma'(t')$, on a
$\underline{\Cal S}_{h(t)}'=\{(\bar\lambda',
\bar\mu')\: \bar\mu'=\Theta'(\bar\lambda',h(t))\}$, o\`u $\bar\nu'=
(\bar\lambda',\bar\mu')$, et l'application de sym\'etrie sera par
d\'efinition la s\'erie formelle $d$-vectorielle ${\Cal
R}_h'(\bar\nu',t):=
\bar\mu'-\sum_{\gamma\in \N^m}\bar{\lambda'}^{\gamma}\Theta_\gamma'
(h(t))\in \C\{\bar\nu'\}\dl t\dr^d$. Dans ces conditions, le lien de
$\Cal R_h'(\bar\nu',t)$ avec $j_{{\bar\nu'}}^\infty
\underline{\Cal S}_{h(t)}'$ s'exprime par une collection infinie
de s\'eries~:
$$
j_{\bar\nu'}^\infty S_{h(t)}'=
(\bar\nu',\{\partial_{\bar\lambda'}^\beta[{\Cal R}_h'(\bar\nu',t)]
\}_{\beta\in\N^m})=(\bar\nu',\{
\partial_{\bar\lambda'}^\beta\left[
\bar\mu'-\Theta'(\bar\lambda',h(t))\right]\}_{\beta\in\N^m}),
\tag 2.2
$$
en prenant les jets {\it avec d\'ependance du point base}, {\it e.g.}
$j_t^k\Psi(t):=(t,\{\partial_t^\alpha\Psi(t)\}_{\v \alpha\v \leq k})$.
Par d\'efinition, cette application formelle ${\Cal R}_h'$ d\'epend du
syst\`eme de coordonn\'ees $t'$, mais sa convergence en est
ind\'ependante (fait qui d\'ecoule de l'invariance biholomorphe des
vari\'et\'es de Segre). Voici notre r\'esultat principal~:

\proclaim{Th\'eor\`eme 2.3} 
Si $(M,p)$ est minimale, alors
l'application de sym\'etrie ${\Cal R}_h'$ est convergente, {\it i.e.}
$\Cal R_h'(\bar\nu',t)\in\C\{\bar\nu',t\}^d$.
\endproclaim

\remark{Remarques} 
{\bf (a)} 
De mani\`ere \'equivalente~: {\it toutes les applications} \,
$\Theta_\gamma'(h(t))=:\varphi_\gamma'(t)$ 
$($une infinit\'e$)$ {\it appartiennent \`a
$\C\{t\}^{d}$, et \, $\exists \ \varepsilon,C>0$, 
$\n t\n \leq \varepsilon \Rightarrow \n \varphi_\gamma'(t)\n\leq
C^{\v \gamma\v+1}$}.

{\bf (b)}
L'\'etude de la r\'egularit\'e de $\Cal R_h'$ {\it g\'en\'eralise
ad\'equatement le principe de sym\'e\-trie de Schwarz en plusieurs
variables}. En effet, aucune condition de non-d\'eg\'en\'eres\-cence
n'est suppos\'ee sur la vari\'et\'e image $(M',p')$, comme dans le
cas $n=1$.

{\bf (c)} 
On supposera dans la suite $(M,p)$ minimale au sens de
Tumanov ({\it voir} [3]).
\endremark

\smallskip
\noindent
{\it Applications.}
Voici maintenant deux applications importantes de ce th\'eor\`eme.

$\bullet$ Premi\`erement, rappelons que $(M',p')$ est
holomorphiquement non-d\'eg\'en\'er\'ee si et seulement si le rang
g\'en\'erique de l'application $t'\mapsto
(\Theta_\gamma'(t'))_{\beta\in\N^m}$ est
\'egal \`a $n$ (crit\`ere de Stanton [19]). Dans ce cas, les
hypoth\`eses du lemme suivant (qui d\'ecoule directement de [1])
seront satisfaites~:

\proclaim{Lemme 2.4} 
Soient $R(t,t')\in \C\{t,t'\}^n$, $t,t'\in \C^n$, et
$h(t)\in\C\dl t\dr^n$, $h(0)=0$, v\'erifiant~: $R(t,h(t)) \equiv 0$ et
$\hbox{\rm d\'et} \, ({\partial R_k\over\partial y_l}(t,h(t)))_{1\leq
k,l\leq n}\not\equiv 0$. Alors $h(t)\in \C\{t\}^n$.
\endproclaim

\noindent
En effet, d'apr\`es [3,19], il existe $\gamma^1,\ldots,
\gamma^n\in\N^m$ tels que $\text{\rm d\'et}
({\partial \Theta_{\gamma^k}'\over
\partial t_l'}(t'))_{1\leq k,l\leq n}\not\equiv 0$, et comme
$\text{\rm d\'et}
({\partial h_k(t)\over \partial t_l})_{1\leq k,l\leq n}(0)\neq 0$,
l'hypoth\`ese est satisfaite avec $R_k(t,t'):=
\Theta_{\gamma^k}'(t')-\varphi_{\gamma^k}'(t)$,
$1\leq k\leq n$, o\`u $\Theta_\gamma'(h(t))=:\varphi_\gamma'(t)\in
\C\{t\}^d$ d'apr\`es le Th\'eor\`eme~2.3.

\smallskip
Par cons\'equent, d'apr\`es le crit\`ere de Stanton et le Lemme~2.4~:

\proclaim{Corollaire 2.5} 
Si $(M',p')$ est holomorphiquement
non-d\'eg\'en\'er\'ee, $h(t)\in\C\{t\}^n$.
\endproclaim

R\'eciproquement, il est connu que si $(M',p')$ est holomorphiquement
d\'eg\'en\'er\'ee, il existe $h^\sharp~: (M',p')\to_{\Cal F}
(M',p')$ inversible et non convergente ({\it cf.} [5] et
\S9 ci-dessous). Le Corollaire~2.5 donne ainsi une {\it 
condition n\'ecessaire et suffisante} pour la convergence de
l'application formelle $h(t)$.

Plus g\'en\'eralement, notre r\'esultat peut aussi
s'interpr\'eter comme un r\'esultat d'analyticit\'e partielle de $h$,
dans l'esprit de Coupet-Pin\-chuk-Sukhov ({\it voir} [8,9])~: 

\proclaim{Corollaire 2.6}
Le degr\'e de transcendance de l'extension de corps
$\hbox{Frac}(\C\{t\})\! \to
\hbox{Frac}(\C(\{t\}))(h_1(t),\ldots,h_n(t)))$ $\subset \hbox{Frac}(\C\dl
t\dr)$ est inf\'erieur ou \'egal au rang g\'en\'erique $e'$ de
l'application holomorphe $t'\mapsto
(\Theta_\gamma'(t'))_{\gamma\in\N^m}$.
\endproclaim
En effet, le graphe formel de $h$ est alors contenu dans une composante
irr\'eductible de dimension $n+e'$ de l'ensemble {\it analytique
complexe} $\{(t,t')\: \Theta_\gamma'(t')=\varphi_\gamma'(t), \,
\forall \, \gamma\}$ et le Corollaire~2.6 d\'ecoule donc de la
caract\'erisation du degr\'e de transcendance par la dimension
minimale d'un ensemble analytique contenant ce graphe
({\it voir} [9]).

\smallskip
$\bullet$
Deuxi\`emement, d'apr\`es le th\'eor\`eme d'approximation d'Artin [1]
appliqu\'e aux \'equations analytiques
$\Theta_\gamma'(h(t))=\varphi_\gamma'(t)\in\C\{t\}^d$ ($\gamma\in
\N^m$), il existe $H(t)\in\C\{t\}^n$ v\'erifiant
$\Theta_\gamma'(H(t))=\varphi_\gamma'(t)\in\C\{t\}^d$, $\forall \
\gamma$. L'application $t\mapsto H(t)$ \'etablissant alors une
\'equivalence convergente ({\it voir} \S8), on en d\'eduit~:

\proclaim{Corollaire 2.7} 
Les vari\'et\'es CR-g\'en\'eriques $\Cal C^\omega$ minimales
$(M,p)$ et $(M',p')$ sont biholomorphes si et seulemement si elles
sont formellement \'equivalentes.
\endproclaim

\remark{Remarque} 
Ce r\'esultat a \'et\'e obtenu r\'ecemment par Baouendi, Rothschild et
Zaitsev dans [6] en supposant l'application $t'\mapsto
(\Theta_\gamma'(t'))_{\gamma\in\N^m}$ de rang constant au voisinage de
$p'$, et en travaillant avec un raffinement du th\'eor\`eme d'Artin
d\^u \`a Wavrick.
\endremark

\subhead 2.8. Conditions de non-d\'eg\'en\'erescence \endsubhead
Ainsi, tout repose sur le Th\'eor\`e\-me 2.3, mais habituellement, les
th\'eor\`emes de r\'egularit\'e d'applications CR sont \'etablis avec
des conditions suppl\'ementai\-res de non-d\'eg\'en\'erescence. Toutes
ces conditions de non-d\'eg\'en\'erescence peuvent s'exprimer \`a
travers le {\it morphisme des $k$-jets de vari\'et\'es de Segre},
introduit par Diederich et Webster ([10]), et qui a servi \`a d\'efinir
l'application de sym\'etrie ci-dessus. Soit $j_{t'}^k{\Cal
S}_{\tau'}'$ le $k$-jet au point $t'$ de la vari\'et\'e de Segre
complexifi\'ee ${\Cal S}_{\tau'}'=\{(w',z')\:
z'=\bar\Theta'(w',\tau')\}$, qui d\'efinit une application holomorphe
$\varphi_k' : {\Cal M}' \ni (t', \tau') \mapsto j_{t'}^k{\Cal
S}_{\tau'}'=(t',\{\partial_{w'}^\beta [z'-\bar\Theta'(w',\tau')]
\}_{\v\beta\v\leq k}) \in \C^{n+d{(m+k)!\over m! \ k!}}$.
Soit ${p'}^c:=(p',\bar p')\in\Cal M'$. Classons toutes ces conditions
par ordre croissant de g\'en\'eralit\'e (pour l'implication {\bf nd4}
$\Rightarrow$ {\bf nd5}, {\it voir} [15]). On dit que $(M',p')$ est~:

\smallskip

{\bf nd1.} {\it Levi-non-d\'eg\'ener\'ee} en $p'$ si
$\varphi_1'$ est une immersion en ${p'}^c$.

{\bf nd2.} {\it Finiment non-d\'eg\'en\'er\'ee} en $p'$ si
$\varphi_k'$ est immersive en ${p'}^c$, $\forall \ k\geq k_0$.

{\bf nd3.} {\it Essentiellement finie} en $p'$ si
$\varphi_k'$ est finie en ${p'}^c$, $\forall \ k\geq k_0$.

{\bf nd4.} {\it Segre-non-d\'eg\'en\'er\'ee} en $p'$ si
$\hbox{rg-g\'en}(\varphi_k'\vert_{\underline{\Cal S}_{p'}'})=m'$, 
$\forall k\geq k_0$.

{\bf nd5.} {\it Holomorphiquement non-d\'eg\'en\'er\'ee} si
$\hbox{rg-g\'en} (\varphi_k')=\dim_\C {\Cal M}'$, 
$\forall \ k\geq k_0$.

\smallskip
\noindent
En particulier, l'\'equivalence formelle $h$ converge sous chacune des
conditions {\bf nd1} \`a {\bf nd5} ci-dessus, d'apr\`es le
Corollaire~2.6, qui g\'en\'eralise des r\'esultats d\'ej\`a connus.

\smallskip
\noindent
\remark{Comparaison avec les r\'esultats pr\'ec\'edents} La 
d\'emonstration du cas {\bf nd1} suit implicitement des formes
normales de Chern-Moser [7]~; le cas {\bf nd2} est trait\'e dans [2]~;
le cas {\bf nd3} dans [4]~; le cas {\bf nd4} par l'auteur dans [15]~;
et enfin, le cas {\bf nd5} par l'auteur, dans [15], en codimension un
(cas hypersurface). Ensuite Mir a trait\'e le Th\'eor\`eme~2.3 en
codimension un ({\it voir} [16]). Mir avait \'ecrit auparavant un
pre\-print traitant le Th\'eor\`eme~2.3 lorsque $(M',p')$ est
alg\'ebrique ({\it cf.}~aussi [9]), cas tr\`es sp\'ecial, puisque
toutes les s\'eries $\{\Theta_\gamma'(t')\}_{\gamma\in\N^m}$ sont {\it
alg\'ebriquement d\'ependantes} sur un nombre {\it fini} d'entre
elles. Cependant, mentionnons que l'utilisation d'{\it identit\'es de
r\'eflexion conjugu\'ees}, qui constitue l'i\-d\'ee principale dans la
d\'emonstration du pr\'esent article ({\it voir} \S5 ci-dessous) et
qui permet de traiter la codimension quelconque sans hypoth\`ese
d'alg\'ebricit\'e, est une id\'ee qui appara\^{\i}t d\'ej\`a
clairement dans [15] sans y avoir \'et\'e toutefois exploit\'ee 
de mani\`ere appropri\'ee.
\endremark

\remark{Remarque sur les conditions de non-d\'eg\'en\'erescence} 
\`A chaque \'etape de la d\'emon\-stration par r\'ecurrence du
Corollaire~2.6, chacune des conditions {\bf nd1}, {\bf nd2}, {\bf nd3}
et {\bf nd4} permet de substituer \`a la consid\'eration de l'{\it
infinit\'e} de s\'eries formelles $\{\Theta_\gamma'(h(t))\}_{\gamma\in\N^m}$
celle des $n$ composantes de $h$. Dans ce cas, il n'est pas n\'ecessaire 
de travailler avec les identit\'es de r\'eflexion conjugu\'ees ({\it
voir} [2,3,4,5,6,16]). Mais ici, dans le cas g\'en\'eral, on
travaillera avec cette collection infinie
$\{\Theta_\gamma'(h(t))\}_{\gamma\in\N^m}$ en utilisant de mani\`ere
cruciale la sym\'etrie du probl\`eme par conjugaison complexe.
\endremark

\smallskip

\remark{Remerciements}
L'auteur tient \`a reconna\^{\i}tre ici sa dette envers les travaux de
Ba\-ouendi-Rothschild et co-auteurs dont il s'est inspir\'e.
\endremark

\head \S3 Notations et pr\'eliminaires \endhead

\subhead 3.1. \'Equations d\'efinissantes \endsubhead
Dans des coordonn\'ees $t=(w,z)\in\C^m\times \C^d$ et
$t'=(w',z')\in\C^m\times \C^d$ s'annulant en $p$ et en $p'$ et telles
que $T_0^cM\cap \C_z^d=\{0\}$ et $T_0^cM'\cap \C_{z'}^d=\{0\}$, 
les vari\'et\'es $(M,0)\subset (\C_t^{n},0)$ et $(M',0)
\subset (\C_{t'}^{n},0)$ sont donn\'ees par deux syst\`emes
de $d$ \'equations scalaires analytiques $z_j=\bar\Theta_j(w,\bar t)$
et $z_j'=\bar\Theta_j'(w',\bar t')$, $1\leq j\leq d$. Soient
$\tau:=(\zeta,\xi):=(\bar t)^c$ et $\tau':=(\zeta',\xi'):=(\bar t')^c$
les variables complexifi\'ees formelles auxquelles correspondent les
complexifications extrins\`eques $(\Cal M,0)\subset \C_t^n\times
\C_\tau^n$ de $(M,0)$ et $(\Cal
M',0)\subset \C_{t'}^n\times \C_{\tau'}^n$ de $(M',0)$, dont voici les
\'equations donn\'ees sous les formes conjugu\'ees {\it \'equivalentes}~:
$$
\Updownarrow
\left\{
\aligned
&
z=\sumg w^\gamma \ \bar\Theta_\gamma(\tau) \ \ \ \ \ 
\hbox{et} \ \ \ \ 
z'=\sumg {w'}^\gamma \ \bar\Theta_\gamma'(\tau'),\\
&
\xi=\sumg \zeta^\gamma \ \Theta_\gamma(t) \ \ \ \ \ 
\hbox{et} \ \ \ \ 
\xi'=\sumg {\zeta'}^\gamma \ \Theta_\gamma'(t').
\endaligned\right.
\tag 3.2
$$
Quitte \`a choisir des coordonn\'ees {\it normales}, on supposera
$\Theta_\gamma(0)=0$ et $\Theta_\gamma'(0)=0$, $\forall \ \gamma$.
Ici, les fonctions analytiques $d$-vectorielles
$\Theta_\gamma(t):=\sum_{\alpha\in\N^n}
\theta_{\gamma,\alpha}t^\alpha\in \C\{t\}^d$, satisfont une in\'egalit\'e 
de Cau\-chy~: $\n \Theta_\gamma(t)\n\leq C^{\v \gamma\v+1}$ pour $\n t\n
\leq \varepsilon$, o\`u $\varepsilon, C>0$, et de m\^eme pour 
$\Theta_\gamma'$. On pose $r(t,\tau):=z-\bar\Theta(w,\tau)$ et
$r'(t',\tau'):=z'-\bar\Theta'(w',\tau')$. Ainsi, $\bar
r(\tau,t)=\xi-\Theta(\zeta,t)$ et $\bar
r'(\tau',t')=\xi'-\Theta'(\zeta',t')$. Par hypoth\`ese, il existe
deux matrices inversibles $a(t,\tau)\in \C\{t,\tau\}^{d\times d}$ et
$a'(t',\tau')\in \C\{t',\tau'\}^{d\times d}$ telles que
$r(t,\tau)\equiv a(t,\tau) \, \bar r(\tau,t)$ et $r'(t',\tau')\equiv
a'(t',\tau') \, \bar r'(\tau',t')$, avec $a(0)=-I_{d\times d}=a'(0)$.

\subhead 3.3. Application formelle \endsubhead
L'application formelle donn\'ee $h(t)$ est par d\'efinition une
s\'erie formelle vectorielle $h(t)=(h_1(t),\ldots,h_n(t))\in\C\dl
t\dr^n$, $h(0)=0$ v\'erifiant $\hbox{d\'et} \, ({\partial h_k\over
\partial t_l})_{1\leq k,l\leq n}(0)\neq 0$ et (apr\`es 
complexification) $r'(h(t),\bar h(\tau))\equiv c(t,\tau) \,
r(t,\tau)$, pour une matrice (n\'ecessairement inversible)
$c(t,\tau)\in\C\dl t,\tau \dr^{d\times d}$. Afin de la distinguer
d'une v\'eritable application ensembliste, on \'ecrira $h^c=(h,\bar
h)\: (\Cal M,0)\to_{\Cal F} (\Cal M',0)$ avec l'indice $\Cal F$ pour
``formel'' ou bien on \'ecrira ``$r'(h(t),\bar h(\tau))=0$ lorsque
$r(t,\tau)=0$''. Une telle expression a un sens, puisque sur la
vari\'et\'e complexe $(2m+d)$-dimensionnelle $(\Cal M,0)=\{(t,\tau)\:
r(t,\tau)=0\}$, on peut choisir (indiff\'eremment) les coordonn\'ees
$(w,\tau)$ ou $(\zeta,t)$ et alors on entendra l'identit\'e
$r'(h(t),\bar h(\zeta,\Theta(\zeta,t)))\equiv 0$ dans $\C\dl
\zeta,t\dr^d$. Bien entendu, on pourrait utiliser le langage
ad\'equat de la th\'eorie des morphismes d'alg\`ebres locales, mais
les calculs pr\'esentant d\'ej\`a une certaine complexit\'e, nous
pr\'ef\'ererons les exposer sous une forme directe. Ainsi,
l'hypoth\`ese $h(\Cal M,0)\subset_\Cal F(\Cal M',0)$ sera exprim\'ee
par les {\it deux \'equations \'equivalentes}~:
$$
f(t)= \bar\Theta'(g(t),\bar h(\tau)) \ \ \ (\Longleftrightarrow)
\ \ \ \bar f(\tau)= \Theta'(\bar g(\tau),h(t)), \ \ \ \ \
\hbox{``lorsque} \ \ r(t,\tau)=0\text{\rm ''}.
\tag 3.4
$$ 
Enfin, nous aurons besoin d'une {\it propri\'et\'e de sym\'etrie par
conjugaison complexe} des objets et des s\'eries formelles. Soit
$h^c=(h,\bar{h})\: (\Cal M,0)\to_{\Cal F} (\Cal M',0)$ comme
ci-dessus, $h^c(t,\tau)=(h(t),\bar{h}(\tau))$ et posons
$\sigma(t,\tau):=(\bar \tau,\bar t)$, $\sigma'(t',\tau'):=(\bar
\tau',\bar t')$. Alors~:
$$
\sigma'(h^c(t,\tau))=\sigma'(h(t),\bar h(\tau))=
(h(\bar \tau), \bar h(\bar t))= h^c(\bar \tau, \bar t)= h^c(\sigma
(t,\tau)).
\tag 3.5
$$

\subhead 3.6.
Application de sym\'etrie \endsubhead
L'{\it application de sym\'etrie} $\Cal R_h' (\bar\nu',t)$, $t\in
\C^n$, $\bar\nu'=(\bar\lambda', \bar\mu')\in \C^{m} \times \C^d$ 
sera par d\'efinition la s\'erie formelle $d$-vectorielle~:
$$
{\Cal R}_h'(\bar\nu',t)=\bar\mu'-\sumg
\bar{\lambda'}^\gamma \ \Theta_\gamma'(h(t))
\in\C\{\bar\nu'\}\dl t\dr^d
\tag 3.7
$$

\remark{Remarque} 
Soient des variables $x_1,x_2\in\C$. L'anneau $\C\dl x_1\dr\{x_2\}$ 
n'a pas de sens.
\endremark

\proclaim{Lemme 3.8}
La propri\'et\'e ${\Cal R}_h'(\bar\nu',t)
\in\C\{\bar\nu',t\}^d$ est ind\'ependante du choix des coordonn\'ees
$(w',z')$ telles que $T_0^cM'\cap \C_{z'}^d=\{0\}$.
\endproclaim

\demo{Preuve}
Cons\'equence ais\'ee de l'invariance biholomorphe des vari\'et\'es de 
Segre.
\qed\enddemo

\head \S4. Convergence de ${\Cal R}_h'$ et de ses jets sur
les cha\^{\i}nes de Segre
\endhead

On note ${\Cal L}:=({\Cal L}^1,\ldots,{\Cal L}^m)$ et
$\underline{\Cal L}:=(\underline{\Cal L}^1,\ldots,\underline{\Cal
L}^m)$ des bases de $T^{1,0}{\Cal M}$ et $T^{0,1}{\Cal M}$ \`a
coefficients holomorphes qui commutent, donn\'ees par $\Cal
L^j:={\partial \over \partial w_j}+{\partial
\bar\Theta (w,\tau)\over \partial w_j} {\partial
\over \partial z}$ et $\underline{\Cal L}^j:={\partial \over \partial
\zeta_j}+{\partial \Theta(\zeta,t)\over \partial \zeta_j}
{\partial \over \partial \xi}$, et ${\Cal L}_w(0)={\Cal L}_{w^1}^1
(\ldots {\Cal L}_{w^m}^m(0))$ le $m$-flot de ${\Cal L}$, $w\in\C^m$
(de m\^eme pour $\underline{\Cal L}_{\zeta}(0)$, $\zeta\in
\C^m$). Clairement, $\Cal L_w(0)=(w,\bar\Theta(w,0),0,0)\in\C^{2n}$ et
$\underline{\Cal L}_\zeta(0)=(0,0,\zeta,\Theta(\zeta,0))\in\C^{2n}$.
Les concat\'enations altern\'ees de tels flots sont
appel\'ees {\it $k$-cha\^{\i}nes de Segre}, par exemple, pour $k=2j$,
$(w_1,\ldots,w_{2j})\mapsto \underline{\Cal L}_{w_{2j}}( {\Cal
L}_{w_{2j-1}}( \ldots \underline{\Cal L}_{w_2}( {\Cal
L}_{w_1}(0))))\in {\Cal M}$, {\it i.e.} $\underline{\Cal L}_{w_2}( {\Cal
L}_{w_1}(0))=(w_1,\bar\Theta(w_1,0),w_2,
\Theta(w_2,$ $w_1,\bar\Theta(w_1,0)))$, {\it etc.}
({\it voir} [14]). On a aussi $\sigma(\Cal
L_w(0))=\underline{\Cal L}_{\bar w}(0)$, {\it etc.}. Si on convient de noter
$w_{(k)}:=(w_1,\ldots,w_k)$, o\`u $w_1,\ldots,w_k\in\C^m$ sont
proches de $0$, ces $k$-cha\^{\i}nes seront abr\'eg\'ees dans la
suite par $\Gamma_k(w_{(k)})$.

\remark{Remarque}
Un point de vue ensembliste \'equivalent sur les cha\^{\i}nes 
de Segre a \'et\'e d\'evelopp\'e par
Baouendi, Ebenfelt et Rothschild ({\it voir} [3]).
Dans [14], l'auteur a ensuite g\'eom\'etris\'e ce point de vue
en introduisant les flots de champs de vecteurs.
\endremark

\smallskip

On note $\nabla_t^\kappa \chi(t):=(\partial_t^\beta \chi(t))_{\v
\beta\v\leq \kappa}\in\C\dl t\dr^{K{(n+\kappa)!\over n! \, \kappa!}}$ 
le $\kappa$-jet d'une~s\'e\-rie formelle vectorielle $\chi(t)\in\C\dl
t\dr^K$, o\`u $K\in\N_*$, $\kappa\in\N$.  Par exemple,
$\nabla_t^\kappa \Cal R_h'(\bar\nu',t)=
\sumg \bar{\lambda'}^\gamma \ \nabla_t^\kappa
[\Theta_\gamma'(h(t))]$.

\proclaim{Lemme 4.1}
Soit $q(x)\in \C\dl x\dr^{2n}$, avec $q(0)=0$. On a $\nabla_t^\kappa
h(\underline{\Cal L}_{\zeta} (q(x)))\equiv
\nabla_t^\kappa h(q(x))$ dans $\C\dl x,\zeta\dr^{n{(n+\kappa)!\over
n! \, \kappa!}}$ et $\nabla_\tau^\kappa \bar{h}({\Cal L}_w 
(q(x))) \equiv \nabla_\tau^\kappa
\bar{h}(q(x))$ dans $\C\dl x,w\dr^{n{(n+\kappa)!\over
n! \, \kappa!}}$.
\endproclaim

\demo{Preuve}
On note $q(x)=(q_1(x),q_2(x))\in\C^n\times \C^n$ et
$q_2(x):=(q_2^1(x),q_2^2(x))\in \C^m\times \C^d$. Il est facile de
voir que $\underline{\Cal L}_\zeta(q(x))=(q_1(x),
\zeta+q_2^1(x),\Theta(\zeta+q_2^1(x),q_1(x)))$, et comme $\nabla_t^\kappa
h(t,\tau)=\nabla_t^\kappa h(t)$, on a bien $\nabla_t^\kappa
h(\underline{\Cal L}_{\zeta} (q(x))) 
\equiv \nabla_t^\kappa h(q_1(x))\equiv
\nabla_t^\kappa h(q(x))$. La deuxi\`eme propri\'et\'e 
d\'ecoule de la premi\`ere gr\^ace \`a la sym\'etrie par conjugaison
complexe~\thetag{3.5} et gr\^ace \`a la relation $\sigma (\Cal
L_w(q(x)))=\underline{\Cal L}_{\bar w}(\sigma(q(x)))$.
\qed
\enddemo

D'apr\`es le crit\`ere de minimalit\'e de [3] reformul\'e dans [14]~:

\proclaim{Lemme 4.2}
La vari\'et\'e CR-g\'en\'erique $\Cal C^\omega$
$(M,p)$ est minimale si et seulement si \, $\Gamma_k\: \C^{mk}\to
(\Cal M,0)$ induit une submersion en $0$ pour $k\geq 2d+1$.
\endproclaim

Ainsi, il nous suffira de d\'emontrer que ${\Cal
R}_h'(\bar\nu',\Gamma_k(w_{(k)}))\in
\C\{\bar\nu',w_{(k)}\}^d$, $\forall \ k\in \N$ pour en d\'eduire 
que $\Cal R_h'(\bar\nu',t)\in\C\{\bar\nu',t\}^d$, comme souhait\'e
({\it cf.}~[2,3,4,6,15,16]).

\subhead 4.3. Estim\'ees de Cauchy \endsubhead
Dans cet objectif, il est n\'ecessaire d'estimer la croissance
des normes $\n \Theta_\gamma'(h(t))\n $ quand $\v \gamma\v\to
\infty$. Commen\c cons par une formule standard et universelle 
de d\'erivation compos\'ee.

\proclaim{Lemme 4.4}
Soit $\beta\in\N^m$.
Il existe un polyn\^ome universel $d$-vectoriel $Q_\beta$ tel que
$$
\left\{
\aligned
&
\partial_t^\beta[\Theta_\gamma'(h(t))]=Q_\beta(\{\partial_t^\delta
h(t)\}_{\delta\leq \beta}, \{[\partial_{t'}^\delta\Theta_\gamma']
(h(t))\}_{\delta\leq \beta})\\
&
Q_\beta(\{\partial_t^\delta
h(t)\}_{\delta\leq \beta},0)\equiv 0.
\endaligned\right.
\tag 4.5
$$
\endproclaim

\demo{Preuve}
On note $Q_\beta(\{h_j^\delta\}_{\delta\leq \beta}^{1\leq j\leq n},
\{{{\theta_k'}^\delta}\}_{\delta\leq \beta}^{1\leq k\leq d})$ 
ce polyn\^ome. Il v\'erifiera \thetag{4.5} si et seulement si
$Q_0={\theta'}^0$ et, par r\'ecurrence,
$$
Q_{\beta+\beta^1}:=\sum_{\delta\leq\beta}\sum_{j=1}^n \
{\partial Q_\beta\over \partial h_j^\delta} \ h_j^{\delta+\beta^1}+
\sum_{\delta\leq \beta} \sum_{k=1}^d \sum_{j=1}^n \
{\partial Q_\beta \over \partial {\theta_k'}^\delta} \
{\theta_k'}^{\delta+\1_j} \ h_j^{\beta^1}
\tag 4.6
$$
o\`u $\beta^1\in\N^m$ avec $\v\beta^1\v=1$ et 
o\`u $\1_j=(0,\ldots,1,\ldots,0)$
avec $1$ \`a la $j$-i\`eme place.
\qed
\enddemo

Si $\Psi(t')\in\C\dl t'\dr^n$, on note $[\nabla_t^\kappa
\Psi(h)](\Gamma_k(w_{(k)})):=
\{[\partial_t^\beta\Psi(h(t))]\v_{t:=\Gamma_k(w_{(k)})}\}_{\v
\beta\v\leq\kappa}$.
Maintenant, en utilisant les polyn\^omes $Q_\beta$, l'estim\'ee de
Cauchy $\n \nabla_{t'}^\kappa
\Theta_\gamma'(t')\n \leq {C'}^{\v\gamma\v+1}$ et le th\'eor\`eme 
d'Artin ([1]), on peut estimer $\n \Theta_\gamma'(h(t))\n $ 
quand $\v \gamma\v\to \infty$~:

\proclaim{Lemme 4.7} 
Soient $k\in\N$ et $\kappa\in\N$. Les propri\'et\'es suivantes
sont \'equivalentes~\text{\rm :}
\roster
\item"{\bf (1)}" \
$[\nabla_t^\kappa{\Cal R}_h'](\bar\nu',\Gamma_k(w_{(k)}))\in
\C\{\bar\nu',w_{(k)}\}^{d{(n+\kappa)!\over
n! \, \kappa!}}$.

\item"{\bf (2)}" \
$[\nabla_t^\kappa
\Theta_\gamma'(h)](\Gamma_k(w_{(k)}))\in\C\{w_{(k)}\}^{d{(n+\kappa)!\over
n! \, \kappa!}}$, $\forall \
\gamma \in\N^m$ et \, $\exists \, C>0$, $\exists \, \varepsilon>0$ tq~{\rm :}

\hskip -0.3cm 
$\n w_{(k)}\n \leq \varepsilon \ \Rightarrow \ \n \, [\nabla_t^\kappa
\Theta_\gamma'(h)](\Gamma_k(w_{(k)})) \, \n \leq C^{\v \gamma \v +1}$
$($estim\'ee de Cauchy$)$.

\item"{\bf (3)}" \
$[\nabla_t^\kappa
\Theta_\gamma'(h)](\Gamma_k(w_{(k)}))
\in\C\{w_{(k)}\}^{d{(n+\kappa)!\over
n! \, \kappa!}}$, $\forall \
\gamma \in\N^m$.
\endroster
\endproclaim

\demo{Preuve}
On a {\bf (1)} $\Longleftrightarrow$ {\bf (2)} et {\bf (2)}
$\Rightarrow$ {\bf (3)}.
\'Etablissons l'implication {\bf (3)} $\Rightarrow$ {\bf (2)}. 
Soit $\beta\in\N^m$, $\v \beta\v \leq \kappa$.  Soient $C>0$ et
$\varepsilon>0$ deux constantes positives dont la valeur pourra varier
suivant le contexte.  Par hypoth\`ese,
$[\partial_t^\beta\Theta_\gamma' (h)](\Gamma_k(w_{(k)}))=:
\varphi_\gamma^\beta(w_{(k)})\in\C\{w_{(k)}\}^d$, $\forall \
\gamma \in\N^m$. D'apr\`es le Lemme~4.4, cette relation peut s'\'ecrire~:
$$
Q_\beta(\{\partial_t^\delta h(\Gamma_k(w_{(k)}))\}_{\delta\leq\beta},
\{[\partial_{t'}^\delta \Theta_\gamma'](h(\Gamma_k(w_{(k)})))\}_{\delta
\leq\beta})=\varphi_\gamma^\beta(w_{(k)}).
\tag 4.8
$$
On l'interpr\`ete comme une relation analytique satisfaite par les
s\'eries formelles $\{\partial_t^\delta
h(\Gamma_k(w_{(k)}))\}_{\delta\leq\beta}$. D'apr\`es le th\'eor\`eme
d'Artin, il existe une solution convergente de l'\'eq.~\thetag{4.8},
soit $\{H^\delta(w_{(k)})\}_{\delta\leq\beta}$, o\`u
$H^\delta(w_{(k)})\in\C\{w_{(k)}\}^n$. Or, gr\^ace \`a la relation
$Q_\beta(\{\partial_t^\delta h(t)\}_{\delta\leq \beta},0)\equiv 0$ et
aux estim\'ees de Cauchy satisfaites par les diff\'erentielles
$\partial_{t'}^\delta \Theta_\gamma'(t')$, $\delta\leq\beta$, on en
d\'eduit ais\'ement une estim\'ee de Cauchy $\n
Q_\beta(\{h_j^\delta\}_{\delta\leq \beta}^{1\leq j\leq n},
\{{{\theta_k'}^\delta}\}_{\delta\leq \beta}^{1\leq k\leq d})\n \leq C^{\v
\gamma\v +1}$ pour $\n (\{h_j^\delta\}_{\delta\leq \beta}^{1\leq j\leq n},
\{{{\theta_k'}^\delta}\}_{\delta\leq \beta}^{1\leq k\leq d})\n\leq 
\varepsilon$. Par composition avec la solution convergente 
$\{H^\delta(w_{(k)})\}_{\delta\leq\beta}$, $H^\delta(0)=0$, 
on en d\'eduit l'estim\'ee de Cauchy souhait\'ee
pour $\n w_{(k)}\n\leq \varepsilon$~:
$$
\n Q_\beta(\{H^\delta(w_{(k)})\}_{\delta\leq\beta},
[\partial_{t'}^\delta\Theta_\gamma']
(H^0(w_{(k)}))\}_{\delta\leq\beta})\n
=\n \varphi_\gamma^\beta(w_{(k)})\n
\leq C^{\v \gamma \v+1}. \qed
\tag 4.9
$$
\enddemo

\noindent
Ainsi, il suffit d'\'etablir {\bf (3)} par r\'ecurrence sur $k$ en
admettant \`a chaque \'etape l'\'equivalence avec {\bf (2)}.
Cette r\'ecurrence proc\`edera en deux moments (\S6 et \S7)~:

\smallskip

\proclaim{\'Etape 1}
Pour tout $k\in\N$, on a~\text{\rm :} 
$$
\left\{
\aligned
&
\left<[\nabla_t^\kappa\Theta_\gamma'(h)](\Gamma_{k}(w_{(k)})) \ 
\text{\rm et} \ 
[\nabla_\tau^\kappa\bar\Theta_\gamma'(\bar h)](\Gamma_{k}(w_{(k)}))
\in\C\{w_{(k)}\}^{d{(n+\kappa)!\over
n! \, \kappa!}}, \forall \ \kappa, \gamma 
\right> \ \ \Longrightarrow \\ 
&
\left<[\Theta_\gamma'(h)](\Gamma_{k+1}(w_{(k+1)})) \ 
\text{\rm et} \ 
[\bar\Theta_\gamma'(\bar
h)](\Gamma_{k+1}(w_{(k+1)}))
\in\C\{w_{(k+1)}\}^{d{(n+\kappa)!\over
n! \, \kappa!}}, \forall \ \gamma\right>.
\endaligned\right.
$$
\endproclaim

\proclaim{\'Etape 2}
En supposant l'\'Etape~1 vraie pour $k-1$,
pour tout $\kappa\in\N$, on a~\text{\rm :} 
$$
\left\{
\aligned
&
\left<[\nabla_t^\kappa\Theta_\gamma'(h)](\Gamma_{k}(w_{(k)})) \ 
\text{\rm et} \ 
[\nabla_\tau^\kappa\bar\Theta_\gamma'(\bar h)](\Gamma_{k}(w_{(k)}))
\in\C\{w_{(k)}\}^{d{(n+\kappa)!\over
n! \, \kappa!}}, \forall \ \gamma\right> \ \ \Longrightarrow \\ 
&
\left<[\nabla_t^{\kappa+1}\Theta_\gamma'(h)](\Gamma_{k}(w_{(k)})) \ 
\text{\rm et} \ 
[\nabla_\tau^{\kappa+1}\bar\Theta_\gamma'(\bar h)](\Gamma_{k}(w_{(k)}))
\in\C\{w_{(k)}\}^{d{(n+\kappa+1)!\over
n! \, (\kappa+1)!}}, \forall \ \gamma \right>.
\endaligned\right.
$$
\endproclaim

\remark{Remarque}
Il est clair que pour $k=0$, les hypoth\`eses de l'\'Etape~1 et 
celles de l'\'Etape~2 sont satisfaites.
\endremark

\head \S5. Sym\'etrie par conjugaison complexe et d\'erivations CR
\endhead

\medskip
\noindent
Mais tout d'abord, si $\beta \in \N^m$, on note
$|\beta|:=|\beta_1|+\cdots+ |\beta_m|$ {\it et} $\underline{\Cal
L}^{\beta}:= {\underline{\Cal L}^1}^{\beta_1} \cdots {\underline{\Cal
L}^m}^{\beta_m}$. Appliquant ces d\'erivations {\it aux deux
identit\'es} $r'(h(t), \bar{h}(\tau))=0$ {\it et} \, $\bar
r'(\bar{h}(\tau),h(t))=0$, $t\in\C^n$, $\tau\in\C^n$, $r(t,\tau)=0$,
on obtient {\it deux familles infinies d'\'equations}, 
deux {\bf identit\'es de r\'eflexion conjugu\'ees}, qui sont satisfaites sur
${\Cal M}$~:
$$
\left\{
\aligned
(*): \ \ f\equiv \bar\Theta'(g,\bar h), \ \ \ \ \ 0 \ \ \equiv & \sumg
g^\gamma \ \underline{\Cal L}^\beta (\bar{\Theta}_\gamma'(\bar h)), \
\ \ \ \ \forall \ \beta \in\N_*^m.\\ (\bar{*}): \ \
\bar{f}=\Theta'(\bar{g},h), \ \ \ \underline{\Cal L}^{\beta} \bar{f} \
\equiv & \sumg \underline{\Cal L}^\beta (\bar{g}^\gamma) \
\Theta_\gamma'(h), \ \ \ \ \ \forall \ \beta \in\N_*^m.
\endaligned
\right.
\tag 5.1
$$
Or, il existe une matrice $a'(t',\tau')\in \C\{ t',\tau'\}^{d\times
d}$, telle que $a'(0,0)=-I_{d\times d}$ et $r'(t',\tau')\equiv
a'(t',\tau') \, \bar r'(\tau',t')$. On en d\'eduit dans 
$\C\{t'\}\dl \tau\dr^d$~:
$$
\left<\underline{\Cal
L}^\beta[ r'(t',\bar h(\tau))]=0, \forall \ \beta \in\N^{m}\right> \
\Longleftrightarrow \ \left<\underline{\Cal L}^\beta [\bar r'(\bar
h(\tau), t')]=0, \forall \ \beta \in \N^{m}\right>,
\tag 5.2
$$
et en particulier, l'\'equivalence des syst\`emes $(*)$ et
$(\bar{*})$, que l'on va exploiter. Mais habituellement, {\it seul
le syst\`eme $(\bar *)$ est consid\'er\'e} 
({\it cf.} [2,3,4,5,6,8,9,10,16,18,19]).

Enfin, on utilisera plusieurs fois deux propri\'et\'es formelles qui
d\'ecoulent de transformations lin\'eaires \'el\'ementaires sur des
syst\`emes trigonaux infinis.

\smallskip

$\bullet$ Premi\`erement, d'apr\`es un calcul
qui utilise l'hypoth\`ese $\hbox{d\'et} \, ({\partial 
h_k\over\partial t_l})_{1\leq k,l\leq n}(0)\neq 0$ et qui est classique
dans les travaux de Baouendi-Rothschild ({\it cf.} [3]), on a~:
$$
\left\{
\aligned
&
\left< 
[\underline{\Cal L}^\beta \bar f ](t,\tau)=\sumg 
[\underline{\Cal L}^\beta
\bar g^\gamma](t,\tau)\ \Theta_\gamma'(t'),
\ \ \ \ \ \forall \ \beta\in \N^m
\right> 
\Longleftrightarrow\\
&
\left< 
\underline{\Omega}_\beta(t,\tau,
\nabla^{\v\beta\v}\bar h(\tau))=
\Theta_\beta'(t')+\sumstg {(\beta+\gamma)! \over
\beta !  \ \gamma !} \ \bar g(\tau)^\gamma \ 
\Theta_{\beta+\gamma}'(t'), \ \ \ \ \
\forall \ \beta\in \N^m,
\right> 
\endaligned\right.
\tag 5.3
$$
o\`u les termes $\underline{\Omega}_\beta$ sont analytiques pr\`es
de $0\times 0\times \nabla^{\v \beta\v} \bar h(0)$ et $r(t,\tau)=0$.

\smallskip

$\bullet$ Deuxi\`emement, on a la r\'esolution formelle
directe du syst\`eme trigonal infini~:
$$
\left\{
\aligned
&
\left<
\psi_\beta+\sumstg {(\beta+\gamma)! \over \beta! \
\gamma!} \ \bar g^\gamma\ \psi_{\beta+\gamma} =
\underline{\omega}_{\beta}, \ \
\forall \ \beta \in \N^m
\right>
\\
&
\Longleftrightarrow 
\left< 
\psi_\beta=\underline{\omega}_\beta+\sumstg
(-1)^{\gamma} {(\beta+\gamma)! \over \beta! \ \gamma!} \ \bar
g^\gamma\ \underline{\omega}_{\beta+\gamma}, \ \ 
\forall \ \beta\in \N^m
\right>.
\endaligned\right.
\tag 5.4
$$
Posons $\underline{\Gamma}_k(w_{(k)}):=\sigma(\Gamma_k(\bar w_{(k)}))$,
{\it i.e.} $\underline{\Gamma}_1(w_1)=\underline{\Cal L}_{w_1}(0)$, 
$\underline{\Gamma}_2(w_{(2)})=\Cal L_{w_2}(\underline{\Cal L}_{w_1}(0))$,
{\it etc.} Alors on a dans 
$\C\{w_{(k)}\}^{d{(n+\kappa)!\over n! \, \kappa !}}$~:
$$
\overline{
[\nabla_t^\kappa \Theta_\gamma'(h)](\Gamma_k(w_{(k)}))}=
\left\{
\aligned
&
[\nabla_\tau^\kappa \bar\Theta_\gamma'(\bar h)](
\underline{\Gamma}_k(\overline{w_{(k)}})), 
\ \ \ \ \text{\rm si} \ k \ \text{\rm
est impair},\\
&
[\nabla_\tau^\kappa \bar\Theta_\gamma'(\bar h)](
\underline{\Gamma}_{k-1}(\overline{w_{(k-1)}})),
\ \ \ \ \text{\rm si} \ k \ \text{\rm
est pair}.
\endaligned\right.
\tag 5.5
$$

\head \S6. Saut \`a la cha\^{\i}ne de Segre sup\'erieure \endhead

\demo{Preuve de l'\'Etape 1}
On traite seulement le cas $k$ pair (le cas $k$ impair est similaire
et s'y ram\`ene formellement gr\^ace \`a la sym\'etrie par conjugaison
complexe \thetag{3-5.5}). Tout d'abord, comme $k$ est pair, on a
d\'ej\`a~: $[\bar\Theta_\gamma'(\bar h)](\Gamma_{k+1}(w_{(k+1)}))= [\bar
\Theta_\gamma'(\bar h)](\Gamma_k(w_{(k)}))]\in\C\{w_{(k)}
\}^{d{(n+\kappa)!\over n! \, \kappa !}}$,
$\forall \ \gamma$, c'est-\`a-dire la deuxi\`eme moiti\'e de la
conclusion de l'\'Etape~1.  Plus g\'en\'eralement, en appliquant
encore le Lemme~4.1, on a~:
$[\nabla_\tau^\kappa\bar\Theta_\gamma'(\bar
h)](\Gamma_{k+1}(w_{(k+1)})) =[\nabla_\tau^\kappa
\bar \Theta_\gamma'(\bar h)](\Gamma_k(w_{(k)}))]\in\C\{w_{(k)}
\}^{d{(n+\kappa)!\over n! \, \kappa !}}$, convergent par 
hypoth\`ese.  Ensuite, les coefficients de $\underline{\Cal L}$
\'etant analytiques, on a facilement~:

\proclaim{Lemme 6.1}
Il existe $P_\beta$ analytique tq. $\underline{\Cal L}^\beta(
\bar\Theta_\gamma'(\bar{h}(\tau)))=P_\beta(t,\tau,[\nabla_\tau^{\v
\beta\v}\bar\Theta_\gamma'(\bar h)])(\tau)$.
\endproclaim

\noindent 
Par cons\'equent, tous les termes suivants sont convergents
($\forall \ \gamma\in\N^m$)~:
$$
\aligned
&
[\underline{\Cal L}^\beta \bar\Theta_\gamma'(\bar h)](
\Gamma_{k+1}(w_{(k+1)})) \! = \!
P_\beta(\Gamma_{k+1}(w_{(k+1)}), \! [\nabla_\tau^{\v \beta\v} 
\bar\Theta_\gamma'(\bar h)](
\Gamma_k(w_{(k)}))) \! \in \! \C\{w_{(k+1)}\}^d.
\endaligned
\tag 6.2
$$
Le th\'eor\`eme d'Artin s'applique donc aux \'equations~$(*)$
\'ecrites au point $(t,\tau):=\Gamma_{k+1}(w_{(k+1)})$, dont les
coefficients sont analytiques gr\^ace \`a~\thetag{6.2}, \'equations qui sont
satisfaites par la s\'erie formelle
$h(\Gamma_{(k+1)}(w_{(k+1)}))$. Ainsi, il existe une solution {\it
convergente} de ces \'eqs.~$(*)$ que l'on note $H(w_{(k+1)})\in
\C\{w_{(k+1)}\}^n$. D'apr\`es l'\'equiva\-lence \thetag{5.2} des
syst\`emes $(*)$ et $(\bar{*})$, cette solution satisfait aussi le
syst\`eme~:
$$
\left\{
\aligned
&
\bar{f}(\Gamma_k(w_{(k)})) \equiv \sumg \bar
g^\gamma(\Gamma_k(w_{(k)}) ) \ \Theta_\gamma'(H(w_{(k+1)})), \\
&
{[}\underline{\Cal L}^\beta\bar{f}](\Gamma_{k+1}(w_{(k+1)})) \equiv \sumg
[\underline{\Cal L}^\beta \bar{g}^\gamma]( \Gamma_{k+1}(w_{(k+1)})) \
\Theta_\gamma'(H(w_{(k+1)})), \ \ \forall \ \beta\in\N_*^m.
\endaligned
\right.
\tag 6.3
$$
Apr\`es la transformation lin\'eaire \thetag{5.3} sur ce syt\`eme et 
sur~$(\bar *)$, on a
$$
\left\{
\aligned
&
\Theta_\beta'(H(w_{(k+1)}))+ \sumstg {(\beta+\gamma)! \over \beta! \
\gamma!} \ \bar g^\gamma(\bar\Gamma_k(w_{(k)})) \
\Theta_{\beta+\gamma}'(H(w_{(k+1)}))=\\ 
&
=\underline{\Omega}_{\beta}(\Gamma_{k+1}(w_{(k+1)}),
\nabla^{\v\beta|}_\tau\bar
h(\bar\Gamma_k(w_{(k)})))= \ \ \ \ \ \ \ \ \ \ \ \ \ \ \ \ \ \ \ \hfill
(\forall \ \beta\in\N^m)\\ 
&
= \Theta_\beta'(h(\Gamma_{k+1}(w_{(k+1)})))+
\sumstg {(\beta+\gamma)! \over \beta! \ \gamma!} \ \bar g^\gamma
(\bar\Gamma_k(w_{(k)})) \
\Theta_{\beta+\gamma}'(h(\Gamma_{k+1}(w_{(k+1)}))).
\endaligned
\right.
\tag 6.4
$$
Pour terminer, on applique aux \'eqs.~\thetag{6.4} la
r\'esolution~\thetag{5.4}. On en d\'eduit que
$\Theta_\beta'(h(\Gamma_{(k+1)}(w_{(k+1)})))\equiv
\Theta_\beta'(H(w_{(k+1)}))\in\C\{w_{(k+1)}\}^d$ converge,
$\forall \ \beta\in \N^m$. C'est la premi\`ere moiti\'e de la conclusion
de l'\'Etape~1.
\qed
\enddemo

\head \S7. R\'ecurrence des jets sur une cha\^{\i}ne de Segre \endhead

\demo{Preuve de l'\'Etape 2}
On traite seulement le cas $k$ impair (le cas $k$ pair est similaire
et s'y ram\`ene par sym\'etrie en utilisant~\thetag{3-5.5}).  Tout
d'abord, comme $k$ est impair, on a d\'ej\`a $[\nabla_\tau^{\kappa+1}
\bar\Theta_\gamma'(\bar h)](\Gamma_{k}(w_{(k)}))= [\nabla_\tau^{\kappa+1}\bar
\Theta_\gamma'(\bar h)](\Gamma_{k-1}(w_{(k-1)}))]\in\C\{w_{(k-1)}
\}^{d{(n+\kappa+1)!\over n! \, (\kappa+1)!}}$
convergent $\forall \ \gamma$, puisque l'on suppose vraie l'\'Etape~1
pour $k-1$. C'est la deuxi\`eme moiti\'e de la conclusion de
l'\'Etape~2. En v\'erit\'e, on va effectuer un calcul direct qui
g\'en\'eralise l'\'Etape~1 du \S6 pour montrer $[\nabla_t^\kappa
\Theta_\gamma'(h)](\Gamma_k(w_{(k)}))\in
\C\{w_{(k)}\}^{d{(n+\kappa)!\over n! \, \kappa!}}$, $\forall \ \kappa$,
en faisant au passage une r\'ecurrence sur les jets ({\it i.e.} sur
$\kappa$) du type de l'\'Etape~2 ({\it voir} \, la preuve du
Lemme~7.20).

Soient donc $\xi\in\C^d$ et $w_{(k)}\in\C^{mk}$. Soit
$\Upsilon$ le $d$-champ de vecteurs tangent \`a $\Cal M$ d\'efini par
$\Upsilon:=\frac{\partial}{\partial z}+{\partial \Theta(\zeta, t)\over
\partial z} \frac{\partial }{\partial \xi}$ et soit
$\underline{\Upsilon}:=\frac{\partial }{\partial \xi}+{\partial
\bar{\Theta}(w, \tau)\over \partial \xi}\frac{\partial }{\partial z}$.
On a $[\Cal L,\underline{\Upsilon}]=0$.

\proclaim{Lemme 7.1}
Les propri\'et\'es suivantes sont \'equivalentes~\text{\rm :}
\roster
\item"{\bf (a)}"
$[\nabla_t^\kappa \Theta_\gamma'(h)](\Gamma_k(w_{(k)}))
\in\C\{w_{(k)}\}^{d{(n+\kappa)!\over n! \, \kappa!}}$, 
$\forall \ \kappa\in\N$.
\item"{\bf (b)}"
$[\Cal L^\delta \underline{\Upsilon}^\alpha 
\Theta_\gamma'(h)](\Gamma_k(w_{(k)}))\in\C\{w_{(k)}\}^d$, 
$\forall \ \delta\in\N^m$, $\forall \ \alpha\in\N^d$.
\item"{\bf (c)}" 
$[\underline{\Upsilon}^\alpha 
\Theta_\gamma'(h)](\Gamma_k(w_{(k)}))\in\C\{w_{(k)}\}^d$, 
$\forall \ \alpha\in\N^d$.
\endroster
\endproclaim

\demo{Preuve}
Comme $\Cal L={\partial\over \partial w}+
{\partial\bar\Theta(w,\tau)\over \partial w} {\partial\over \partial
z}$, $\underline{\Upsilon}=\frac{\partial }{\partial \xi}+ {\partial
\bar{\Theta}(w, \tau)\over \partial \xi}\frac{\partial }{\partial z}$,
$\partial_t=(\partial_w,\partial_z)$ et ${\partial
\bar{\Theta}(0,0)\over \partial \xi}=I_{d\times d}$, on a
{\bf (a)} $\Longleftrightarrow$ {\bf (b)} par transformations
lin\'eaires. De plus, on a {\bf (c)} $\Rightarrow$ {\bf (b)},
car $\partial_{w_k}^\delta [\Psi(\Gamma_k(w_{(k)}))]=
\partial_{w_k}^\delta
[\Psi(\Cal L_{w_k}(\Gamma_{k-1}(w_{(k-1)})))]= [\Cal L^\delta
\Psi](\Gamma_k(w_{(k)}))$.
\qed
\enddemo

Il suffit donc de prouver {\bf (c)}.
Soit $(\xi,p)\mapsto \underline{\Upsilon}_\xi(p)$ le $d$-flot de
$\underline{\Upsilon}$. Bien s\^ur, on a
$\partial_\xi^\alpha[\Psi(\underline{\Upsilon}_\xi(q(x)))]
=[\underline{\Upsilon}^\alpha\Psi](\underline{\Upsilon}_\xi(q(x)))$. 
Consid\'erons $\underline{\Upsilon}_\xi(\Gamma_k(w_{(k)}))\in\Cal M$. On pose~:
$$
\left\{
\aligned
&
E_\beta(w_{(k)},\xi,t'):=\sumg
{w'}^\gamma \ [\underline{\Cal L}^\beta \bar\Theta_\gamma'(
\bar h)](\underline{\Upsilon}_\xi(
\Gamma_k(w_{(k)}))),\\
&
F_\beta(w_{(k)},\xi,t'):=
[\underline{\Cal L}^\beta \bar f](\underline{\Upsilon}_\xi(
\Gamma_k(w_{(k)})))-\sumg
[\underline{\Cal L}^\beta \bar g^\gamma](\underline{\Upsilon}_\xi(
\Gamma_k(w_{(k)})))\ 
\Theta_\gamma'(t'), 
\endaligned\right.
\tag 7.2
$$
pour tout $\beta\in\N^m$ ({\it cf.}  $(*)$-$(\bar *)$). D'apr\`es la
d\'efinition de $a'(t',\tau')$ ({\it voir} \S3.1), on a~:
$$
\left\{
\aligned
&
E_0(w_{(k)},\xi,t')=a'(t',\bar h(\underline{\Upsilon}_\xi(
\Gamma_k(w_{(k)})))) \
F_0(w_{(k)},\xi,t'),\\
&
F_0(w_{(k)},\xi,t')=\bar a'(\bar h(\underline{\Upsilon}_\xi(
\Gamma_k(w_{(k)}))),t') \
E_0(w_{(k)},\xi,t').
\endaligned\right.
\tag 7.3
$$
Appliquant toutes les d\'erivations $\underline{\Cal L}^\beta$ aux
\'eqs.~\thetag{7.3}~:
$$
\left\{
\aligned
&
E_\beta(w_{(k)},\xi,t')\equiv a'(t',w_{(k)},\xi) \
F_\beta(w_{(k)},\xi,t')
+\sum_{\delta<\beta}
{a_\delta'}^\beta(t',w_{(k)},\xi) \ F_\delta(w_{(k)},\xi,t'),\\
&
F_\beta(w_{(k)},\xi,t')\equiv b'(t',w_{(k)},\xi) \
E_\beta(w_{(k)},\xi,t')+\sum_{\delta<\beta}
{b_\delta'}^\beta(t',w_{(k)},\xi) \ E_\delta(w_{(k)},\xi,t').
\endaligned\right.
\tag 7.4
$$
Ici, ${a_\delta'}^\beta(t',w_{(k)},\xi)$, ${b_\delta'}^\beta(t',w_{(k)},\xi)
\in\C\{t'\}\dl w_{(k)},\xi\dr^{d\times d}$. 
D'apr\`es~\thetag{7.4}, on a ({\it cf.}~\thetag{5.2}) $\text{\rm
Id\'eal}\left<\{E_\beta(w_{(k)},\xi,t')\}_{\beta\in\N^m} \right>=
\text{\rm Id\'eal}\left<\{F_\beta(w_{(k)},\xi,t')\}_{\beta\in\N^m} \right>$.
Plus g\'en\'eralement, par r\'ecurrence sur $\alpha\in\N^d$, on
d\'efinit une collection $\{E_\beta^{(\alpha)}\}_{\alpha\in\N^d,
\beta\in\N^m}$ de fonctions $d$-vectorielles comme suit. Soit
$\alpha^1\in\N^d$ avec $\v \alpha^1\v =1$. On pose $T_0'=t'$ et~:
$$
\left\{
\aligned
&
E_\beta^{(0)}(w_{(k)},\xi,t'):=
E_\beta(w_{(k)},\xi,t'); \ \ \text{\rm et~:} \
E_\beta^{(\alpha+\alpha^1)}(w_{(k)},\xi,
\{T_{\alpha'}'\}_{\alpha'\leq \alpha+\alpha^1}):=\\
&
:={\partial E_\beta^{(\alpha)}\over \partial \xi^{\alpha^1}}
(w_{(k)},\xi,\{T_{\alpha'}'\}_{\alpha'\leq \alpha})+
\sum_{\alpha'\leq\alpha} {\partial E_{\beta}^{(\alpha)}\over
\partial T_{\alpha'}'}
(w_{(k)},\xi,\{T_{\alpha'}'\}_{\alpha'\leq \alpha}) \
T_{\alpha'+\alpha^1}'.
\endaligned \right.
\tag 7.5
$$
On d\'efinit aussi la collection similaire
$\{F_\beta^{(\alpha)}\}_{\alpha\in\N^d, \beta\in\N^m}$.
Par construction~:
$$
\left\{
\aligned
&
\left[E_\beta^{(\alpha)}(w_{(k)},\xi,\{T_{\alpha'}'\}_{\alpha' \leq 
\alpha})\right]\left\vert_{T_{\alpha'}':=
[\underline{\Upsilon}_\xi^{\alpha'} h](
\underline{\Upsilon}_\xi(\Gamma_k(w_{(k)}))), \
\forall \ \alpha'\leq \alpha}=\right.\\
&
=\partial_\xi^\alpha[E_\beta(w_{(k)},\xi,h(
\underline{\Upsilon}_\xi(\Gamma_k(w_{(k)}))))].
\endaligned\right.
\tag 7.6
$$
Voici maintenant une propri\'et\'e g\'en\'eralisant~\thetag{7.4} qui
se v\'erifie par un calcul formel direct en utilisant les
relations~\thetag{7.3} et les d\'efinitions~\thetag{7.5} de
$E_{\beta}^{(\alpha)}$ et de $F_{\beta}^{(\alpha)}$~: \qed

\proclaim{Lemme 7.7}
Dans l'anneau \, $\C\{ \, \! \{T_{\alpha'}'\}_{\alpha'\leq \alpha} \}
\dl w_{(k)},\xi\dr^d$,
on a pour tout $\alpha\in\N^d$~:
$$
\text{\rm Id\'eal} \left< \{E_\beta^{(\alpha)}(w_{(k)},\xi,
\{T_{\alpha'}'\}_{\alpha'\leq \alpha})\}_{\beta\in\N^m}\right> \! =
\text{\rm Id\'eal} \left< \{F_\beta^{(\alpha)}(w_{(k)},\xi,
\{T_{\alpha'}'\}_{\alpha'\leq \alpha})\}_{\beta\in\N^m}\right>.
\tag 7.8
$$
\endproclaim

\noindent
{\it Suite de la d\'emonstration.}  D'apr\`es $(*)$, la s\'erie
formelle $h(\underline{\Upsilon}_\xi(\Gamma_k(w_{(k)})))$ est une
solution des \'equations $E_\beta(w_{(k)},\xi,$ $h(
\underline{\Upsilon}_\xi(\Gamma_k(w_{(k)}))
))\equiv 0$, $\forall \ \beta\in\N^m$. Par cons\'equent~:
$$
\left\{
\aligned
&
0\equiv\partial_\xi^\alpha \v_{\xi=0} [E_\beta(w_{(k)},\xi,h(
\underline{\Upsilon}_\xi(\Gamma_k(w_{(k)}))))]=\\
&
=E_\beta^{(\alpha)}(w_{(k)},0,
\{[\underline{\Upsilon}^{\alpha'} h]
(\Gamma_k(w_{(k)}))\}_{\alpha'
\leq \alpha}), \ \ \ \ \ 
\forall \ \alpha\in\N^d, \ \forall \ \beta\in\N^m.
\endaligned\right.
\tag 7.9
$$
Maintenant, on fixe $\alpha\in\N^d$ et on consid\`ere 
le sous-syst\`eme fini~:
$$
E_\beta^{(\alpha')}(w_{(k)},0,
\{[\underline{\Upsilon}^{\alpha''} h]
(\Gamma_k(w_{(k)}))\}_{\alpha''
\leq \alpha'})=0, \ \ \ \ \ \forall \ \alpha'\leq \alpha.
\tag 7.10
$$

\proclaim{Lemme 7.11}
Les \'equations~\thetag{7.10} sont analytiques, {\it i.e.}~: 
$$
E_\beta^{(\alpha')}(w_{(k)},0,\{T_{\alpha''}\}_{\alpha''
\leq \alpha'})\in\C\{w_{(k)},\{T_{\alpha''}\}_{\alpha''
\leq \alpha'}\}^d, \ \ \ \ \ \forall \ \alpha'\leq \alpha,
\tag 7.12
$$
\endproclaim

\demo{Preuve}
En effet, comme $k$ est impair, on a $\Gamma_k(w_{(k)})={\Cal
L}_{w_k}(\Gamma_{k-1}(w_{(k-1)}))$, d'o\`u $[\nabla_\tau^\kappa
\bar\Theta_\gamma'(\bar h)](\Gamma_k(w_{(k)}))\equiv
[\nabla_\tau^\kappa \bar\Theta_\gamma'(\bar h)](\Gamma_{k-1}
(w_{(k-1)}))\in\C\{w_{(k-1)}\}^{d{(n+\kappa)!\over n! \, \kappa !}}$,
puisqu'on suppose vraie l'\'Etape~1 pour $k-1$.  Par cons\'equent, en
appliquant le Lemme~6.1, on voit que les d\'eriv\'ees
$\partial_\xi^{\alpha''}\v_{\xi=0}[[\underline{\Cal L}^\beta
\bar\Theta_\gamma'(\bar h)](\underline{\Upsilon}_\xi(\Gamma_k(w_{(k)})
))]\in\C\{w_{(k)}\}^d$ des coefficients de $E_\beta^{(\alpha')}$
({\it cf.} \'eq.~\thetag{7.2}${}^{1^{i\grave{e}re}}$) par 
rapport \`a $\{T_{\alpha''}\}_{\alpha''\leq \alpha}$ convergent toutes.
\qed
\enddemo

Ainsi, il existe des solutions
$H_{\alpha'}(w_{(k)})\in\C\{w_{(k)}\}^n$, $\alpha'\leq \alpha$,
satisfaisant~:
$$
E_\beta^{(\alpha')}(w_{(k)},0,
\{H_{\alpha''}(w_{(k)})\}_{\alpha''\leq \alpha'})\equiv 0,
\ \ \ \ \ \forall \ \alpha'\leq \alpha.
\tag 7.13
$$
Gr\^ace \`a la propri\'et\'e~\thetag{7.8}, on d\'eduit
de~\thetag{7.13}~:
$$
F_\beta^{(\alpha')}(w_{(k)},0,\{H_{\alpha''}(w_{(k)})\}_{\alpha''
\leq \alpha'})
\equiv 0, \ \ \ \ \ \forall \ \alpha'\leq \alpha.
\tag 7.14
$$
Mais on a aussi d'un autre c\^ot\'e en 
d\'erivant~\thetag{7.2}${}^{2^{i\grave{e}me}}$ 
par rapport \`a $\xi$~:
$$
F_\beta^{(\alpha')}(w_{(k)},0,
\{[\underline{\Upsilon}^{\alpha''}h](\Gamma_k(w_{(k)}))\}_{\alpha''
\leq \alpha'})\equiv 0, \ \ \ \ \ \forall \ \alpha'\leq \alpha,
\tag 7.15
$$
une identit\'e que l'on peut r\'e\'ecrire plus explicitement comme suit~:
$$
\left\{
\aligned
&
\partial_\xi^{\alpha'}\v_{\xi=0}
[[\underline{\Cal L}^\beta \bar f](
\underline{\Upsilon}_\xi(\Gamma_k(w_{(k)})))]\equiv \sumg
\sum_{\alpha''\leq\alpha'} 
{\alpha'!\over \alpha''! \ (\alpha'-\alpha'')!} \\
&
\left.
\partial_\xi^{\alpha'-\alpha''}[[\underline{\Cal L}^\beta
\bar g^\gamma](
\underline{\Upsilon}_\xi(\Gamma_k(w_{(k)})))]
\ \partial_\xi^{\alpha''}[\Theta_\gamma'(h(
\underline{\Upsilon}_\xi(\Gamma_k(w_{(k)}))
))]\right\vert_{\xi=0}, \ \ \forall \ \alpha'\leq \alpha.
\endaligned\right.
\tag 7.16
$$ 
Or, il existe clairement des fonctions ${\Theta_\gamma'}^{\alpha''}$
analytiques telles que~:
$$
{\Theta_\gamma'}^{\alpha''}(\{[\underline{\Upsilon}^{\alpha'''}
h](\Gamma_k(w_{(k)}))\}_{\alpha'''\leq \alpha''}):=
\partial_\xi^{\alpha''}\v_{\xi=0} [\Theta_\gamma'(h(
\underline{\Upsilon}_\xi(\Gamma_k(w_{(k)}))))].
\tag 7.17
$$
Par cons\'equent, on peut r\'e\'ecrire les 
\'equations~\thetag{7.14} et \thetag{7.16} comme suit~:
$$
\left\{
\aligned
&
\partial_\xi^{\alpha'}\v_{\xi=0}
[[\underline{\Cal L}^\beta \bar f](
\underline{\Upsilon}_\xi(\Gamma_k(w_{(k)})))]\equiv\sumg
\sum_{\alpha''\leq\alpha'} 
{\alpha'!\over \alpha''! \ (\alpha'-\alpha'')!} \\
&
\partial_\xi^{\alpha'-\alpha''}
\v_{\xi=0}[[\underline{\Cal L}^\beta
\bar g^\gamma](
\underline{\Upsilon}_\xi(\Gamma_k(w_{(k)})))]
\ {\Theta_\gamma'}^{\alpha''}(\{
H_{\alpha'''}(w_{(k)})\}_{\alpha'''\leq \alpha''}), 
\ \ \forall \ \alpha'\leq \alpha.
\endaligned \right.
\tag 7.18
$$
$$
\left\{
\aligned
&
\partial_\xi^{\alpha'}\v_{\xi=0}
[[\underline{\Cal L}^\beta \bar f](
\underline{\Upsilon}_\xi(\Gamma_k(w_{(k)})))]\equiv\sumg
\sum_{\alpha''\leq\alpha'} 
{\alpha'!\over \alpha''! \ (\alpha'-\alpha'')!} \\
&
\partial_\xi^{\alpha'-\alpha''}\v_{\xi=0}[[\underline{\Cal L}^\beta
\bar g^\gamma](
\underline{\Upsilon}_\xi(\Gamma_k(w_{(k)})))]
\ {\Theta_\gamma'}^{\alpha''}(\{[\underline{\Upsilon}^{\alpha'''}
h](\Gamma_k(w_{(k)}))\}_{\alpha'''\leq \alpha''})
\equiv 0.
\endaligned \right.
\tag 7.19
$$
On va maintenant utiliser \thetag{5.3-4} afin de d\'eduire des
\'equations~\thetag{7.18-19}~:

\proclaim{Lemme 7.20}
Pour tout $\gamma\in\N^m$ et tout $\alpha'\leq\alpha$, on a~:
$$
{\Theta_\gamma'}^{\alpha'}(\{[\underline{\Upsilon}^{\alpha''}
h](\Gamma_k(w_{(k)}))\}_{\alpha''\leq \alpha'})\equiv
{\Theta_\gamma'}^{\alpha'}(\{H_{\alpha''}(w_{(k)})\}_{\alpha''\leq \alpha'})
\in\C\{w_{(k)}\}^d.
\tag 7.21
$$
\endproclaim
\demo{Preuve}
En appliquant directement~\thetag{5.3-4}, ceci est vrai pour
$\v\alpha'\v=0$ en consid\'erant les
\'eqs.~\thetag{7.18-19} seulement au rang $\alpha'=0$,
comme \`a la fin du \S6.  Supposons par r\'ecurrence que~\thetag{7.21}
est vrai pour $\alpha'<\alpha$, $\v\alpha'\v=\kappa\in\N_*$. 
Soit $\alpha_0'\leq \alpha$, $\v\alpha_0'\v=\kappa+1$. On \'ecrit les
\'eqs.~\thetag{7.18-19} au rang 
$\alpha':=\alpha_0'$ et on les soustrait deux \`a deux. Gr\^ace \`a
cette hypoth\`ese de r\'ecurrence, on obtient~:
$$
\left\{
\aligned
&
\sumg [\underline{\Cal L}^\beta\bar g^\gamma](\Gamma_{(k)}(w_{(k)})) \
\left[
{\Theta_\gamma'}^{\alpha_0'}(\{
H_{\alpha''}(w_{(k)})\}_{\alpha''\leq \alpha_0'})-\right.\\
&
\left.-{\Theta_\gamma'}^{\alpha_0'}(\{
[\underline{\Upsilon}^{\alpha''}
h](\Gamma_k(w_{(k)}))\}_{\alpha''\leq \alpha_0'})
\right]\equiv 0,
\endaligned\right.
\tag 7.22
$$
pour tout $\beta\in\N^m$. D'apr\`es~\thetag{5.3-4}, on a
alors~\thetag{7.21} pour $\alpha'=\alpha_0'$. Ainsi,
l'\'eq.~\thetag{7.21} conclut que $[\underline{\Upsilon}^\alpha
\Theta_\gamma'(h)](\Gamma_k(w_{(k)}))\in\C\{w_{(k)}\}^d$, 
$\forall \ \alpha\in\N^d$.
\qed
\enddemo

\proclaim{Conclusion 7.23} 
Gr\^ace au Lemme~7.1, on a ainsi achev\'e d'\'etablir la premi\`ere
(et la seconde) moiti\'e de l'\'Etape~2 (cas $k$ impair).  En
appliquant le Lemme~4.7 pas \`a pas dans le proc\'ed\'e de
r\'ecurrence en deux moments d\'efini par l'\'Etape~1 et par
l'\'Etape~2, on en conclut que ${\Cal
R}_h'(\bar\nu',\Gamma_k(w_{(k)}))\in
\C\{\bar\nu',w_{(k)}\}^d$, $\forall \ k\in \N$, et donc 
$\Cal R_h'(\bar\nu',t)\in\C\{\bar\nu',t\}^d$ par minimalit\'e de
$(\Cal M,0)$.
\endproclaim
La d\'emonstration du Th\'eor\`eme~2.3 est termin\'ee.
\qed
\enddemo

\head \S8. \'Equivalences formelles et \'equivalences holomorphes \endhead

L'\'enonc\'e suivant pr\'ecise le contenu du Corollaire~2.7~:

\proclaim{Th\'eor\`eme 8.1}
Soit $h\: (M,0)\to_{\Cal F}(M',0)$ une \'equivalence formelle entre
sous-vari\'et\'es CR-g\'en\'eriques minimales $\Cal C^\omega$ de
$\C^n$. Pour tout $N\in\N_*$, il existe une \'equivalence holomorphe
$H_N\: (M,0)\to_{\Cal F}(M',0)$ avec $H_N(t)\equiv h(t) \ (\hbox{mod}
\ \n t\n^N)$.
\endproclaim

\demo{Preuve}
D'apr\`es le Th\'eor\`eme~2.3, il existe
$\varphi_\gamma'(t)\in\C\{t\}^d$ tel que $\Theta_\gamma'(h(t))\equiv
\varphi_\gamma'(t)$, $\forall \
\gamma\in\N^m$. Soit $N\in\N_*$. Appliquant le th\'eor\`eme d'Artin,
on obtient une application holomorphe $H_N$ satisfaisant
$\Theta_\gamma'(H_N(t))\equiv \varphi_\gamma'(t)$, $\forall \
\gamma\in\N_m$ et $H_N(t)\equiv h(t) \ (\hbox{mod}
\ \n t\n^N)$. Bien s\^ur, on a l'estim\'ee de Cauchy
$\n \Theta_\gamma'(H_N(t))\n \leq C^{\v \gamma\v+1}$ et $\bar
f(\zeta,\Theta(\zeta,t))\equiv
\sumg\bar g^\gamma(\zeta,\Theta(\zeta,t)) \ \Theta_\gamma'(H_N(t))$.
Rappelons $-r'(H_N(t),\bar h(\tau))\equiv a'(H_N(t),\bar h(\tau))
\ \bar r'(\bar h(\tau),H_N(t))$, d'o\`u $F_N(t)\equiv
\sumg G_N^\gamma(t) \ \bar\Theta_\gamma'(
\bar h(\zeta,\Theta(\zeta,t)))$, qui \'equivaut encore \`a
$F_N(w,\bar\Theta(w,\tau)) )\equiv \sumg
G_N^\gamma(w,\bar\Theta(w,\tau)) \
\bar\Theta_\gamma'(\bar h(\tau))$, et finalement 
$F_N(w,\bar\Theta(w,\tau))\equiv
\sumg G_N^\gamma(w,\bar\Theta(w,\tau)) \ \bar\Theta_\gamma'(\bar
H_N(\tau))$. En conclusion, l'ap\-plication $(H_N,\bar H_N)\: (\Cal
M,0)\to (\Cal M',0)$ \'etablit une \'equivalence convergente.
\qed
\enddemo

\head \S9 N\'ecessit\'e de la non-d\'eg\'en\'erescence holomorphe 
\endhead

On rappelle qu'une hypersurface $\Cal C^\omega$ $(M',p')$ est {\it
holomorphiquement d\'eg\'en\'er\'ee} $p'$ si et seulement si il existe
un germe de champ de vecteurs $L'=\sum_{j=1}^n a_j'(t'){\partial\over
\partial t_j'}$ \`a coefficients holomorphes non
tous nuls, {\it tangent \`a $(M',p')$} ({\it voir} [3,19]). La
n\'ecessit\'e de la non-d\'eg\'en\'erescence holomorphe pour la
r\'egularit\'e de $h$ a \'et\'e \'etablie dans [4] pour les
diff\'eomorphismes CR $\Cal C^\infty$, mais l'auteur ne conna\^{\i}t
pas de r\'ef\'erence publi\'ee pour la n\'ecessit\'e dans le cas
formel. En voici une d\'emonstration br\`eve utilisant [1].

\proclaim{Proposition 9.1}
Il existe $\varpi'(t')\in \C\dl t'\dr\backslash \C\{t'\}$,
$\varpi'(0)=0$, tel que le flot $\C^n\ni t'\mapsto \exp (\varpi' (t')
L')(t')\in \C^n$ induit une auto-application formelle inversible
non-convergente $h^\sharp\: (M',0)\to_{\Cal F} (M',0)$.
\endproclaim

\demo{Preuve}
Soit $\varphi'\: (t',u')\mapsto \exp (u'L')(t')=\varphi'(t',u')$, le
flot local de $L'$, qui est holomorphe en $t'\in \C^n$ et $u'\in \C$,
pour $\n t'\n$, $\vert u'\vert\leq \varepsilon$, $\varepsilon>0$. Ce
flot satisfait $\varphi'(t',0)\equiv t'$ et
$\partial_{u'}\varphi_k'(t',u')\equiv a_k'(\varphi'(t',u'))$. Comme
$L'\neq 0$, on a $\partial_{u'} \varphi'(t',u')\not\equiv 0$. On peut
supposer $\partial_{u'}\varphi_1'(t',u')\not\equiv 0$. Soit
$\varpi'(t')\in\C\dl t'\dr\backslash \C\{t'\}$, $\varpi'(0)=0$, une
s\'erie formelle {\it non convergente}, satisfaisant de plus
$\partial_{u'} \varphi_1'(t',\varpi'(t'))\not\equiv 0$ dans $\C\dl t'
\dr$ (il en existe beaucoup). Si la s\'erie formelle $h^\sharp\: 
t'\mapsto_{\Cal
F}\varphi'(t',\varpi'(t'))$ \'etait convergente, alors
$t'\mapsto_{\Cal F} \varpi'(t')$ le serait aussi (par le Lemme~2.4),
contrairement au choix de $\varpi'$. Enfin, $L'$ \'etant tangent 
\`a $(M',0)$, il est clair que $h^\sharp(M',0)\subset_{\Cal F}
(M',0)$.
\qed
\enddemo

\footnotetext[1]{pdf file~: 
cmi.univ-mrs.fr/$\sim$merker/index.html.}

\enddocument
\end